\documentclass{amsart}

\usepackage{amsmath,amsthm,amssymb,amscd,enumerate,mathtools,enumitem}
\usepackage{amsfonts}
\usepackage{rotating}
\usepackage{euscript}
\usepackage{pst-node}
\usepackage{epsfig,verbatim}

\def\be{\begin{equation}}
\def\ee{\end{equation}}

\def\f{\EuScript}
\def\N{{\mathbb N}}

\def\ord{{\rm ord\,}}

\def\phi{{\varphi}}
\def\v{{\varepsilon}}

\def\GCD{{\rm GCD }}

\def\mod{{\rm mod\ }}

\def\bp{\begin{proposition}}
\def\ep{\end{proposition}}

\def\bt{\begin{theorem}}
\def\et{\end{theorem}}
\def\br{\begin{remark}}
\def\er{\end{remark}}
\def\be{\begin{equation}}
\def\bee{\begin{equation*}}
\def\l{\label}
\def\la{\label}

\def\ee{\end{equation}}
\def\eee{\end{equation*}}
\def\bl{\begin{lemma}}
\def\el{\end{lemma}}
\def\bc{\begin{corollary}}
\def\ec{\end{corollary}}
\def\pr{\noindent{\it Proof. }}

\def\bd{\begin{definition}}
\def\ed{\end{definition}}

\newtheorem{theorem}{Theorem}[section]
\newtheorem{lemma}[theorem]{Lemma}
\newtheorem{definition}[theorem]{Definition}
\newtheorem{corollary}[theorem]{Corollary}
\newtheorem{proposition}[theorem]{Proposition}
\newtheorem{problem}[theorem]{Problem}

\theoremstyle{definition}

\theoremstyle{definition}
\newtheorem{remark}[theorem]{Remark}
\def\bpr{\begin{problem}}
\def\epr{\end{problem}}

\mathtoolsset{showonlyrefs}

\begin{document}
\title{Functional equations in formal power series}
\author{Fedor Pakovich}
\thanks{
This research was supported by ISF Grant No. 1092/22}
\address{Department of Mathematics, Ben Gurion University of the Negev, Israel}
\email{
pakovich@math.bgu.ac.il}

\begin{abstract}
Let $k$ be an algebraically closed field of characteristic zero, and $k[[z]]$ the  ring of formal power series over $k$. In this paper, we study 
equations in the semigroup  $z^2k[[z]]$ with the semigroup operation being composition. We prove a number of general results about such equations and provide some applications. 
In particular, we answer a question of Horwitz and Rubel 	about decompositions of ``even'' formal power series.  We also  
show that every right amenable 
subsemigroup of $z^2k[[z]]$ is conjugate to a subsemigroup of the semigroup of monomials.

\end{abstract}

\maketitle

\section{Introduction} 

Let $k$ be an algebraically closed field of characteristic zero, and $k[[z]]$ the ring of formal power series over $k$. For an element $A(z)=\sum_{n\geq 0}c_nz^n$ of $k[[z]]$, we define its {\it order} 
by the formula $\ord A=\min\{n\geq 0\,\vert\, c_n\neq 0\}.$ We denote by $k_n[[z]]$, $n\geq 0$, the subset of  $k[[z]]$ consisting of formal power series of order $n$, and by $\Gamma$ the  subset of  $k[[z]]$ 
consisting of formal power series of order at least two. 
If $A$ and $B$ are elements of $k[[z]]$ with $\ord B\geq 1$, then the operation $A\circ B$ of {\it composition} of $A$ and $B$ is well defined. In particular,  with respect to this operation,    
the set $k_1[[z]]$  is a {\it group}, and  the set $\Gamma$ is a {\it semigroup}.

The group  $k_1[[z]]$ has been intensively studied (see e. g.  \cite{bab2}, \cite{bab1}, \cite{bak}, \cite{bru}, \cite{j}, \cite{1},   \cite{2}, \cite{muc}, \cite{rei}, \cite{sche}, \cite{schw}). In this paper, we focus on the less  studied semigroup $\Gamma$ with an emphasis on equations  in $\Gamma$. In other words, we study  functional equations in formal powers series of order at least two. 
An example of such an equation is simply the equation \be \l{abo} A=A_1\circ A_2 \circ \dots \circ A_r, \ \ \ r\geq 2,\ee where $A\in \Gamma$ is a given and $A_1, A_2, \dots, A_r\in \Gamma$ are unknown, describing the ways in which an element $A$ of $\Gamma$ can be represented as a composition of other elements of  $\Gamma$. 
Although the problem of characterizing solutions of \eqref{abo} is fundamental, 
we were unable to find relevant references in the literature, and provide an answer in this paper. 
Specifically, we describe {\it equivalence classes} of decompositions \eqref{abo}, where two decompositions 
\be \la{decc} A=A_1\circ A_2\circ \dots \circ A_k \ \ \ \ {\rm and} \ \ \ \
A=\widehat A_1\circ \widehat A_{2}\circ \dots \circ \widehat A_m,\ee  
are considered as equivalent if  $k=m$ and there exist elements 
$\mu_i,$ $1\leq i \leq k-1,$ of $k_1[[z]]$ such that
\be \la{decc0} A_1=\widehat A_1\circ \mu_1^{-1}, \ \ \
A_i=\mu_{i-1}\circ \widehat A_i \circ \mu_{i}^{-1}, \ \ \ 1<i< k, \ \ \ {\rm and} \ \ \ A_k=\mu_{k-1}\circ \widehat A_k.
\ee

Let us recall that for every $A\in\Gamma$ of order $n$  there exists an element $\beta_A$ of $ k_1[[z]]$, called the {\it B\"ottcher function}, 
such that 
$$ \beta_A^{-1}\circ  A\circ \beta_A=z^n.$$ 
The B\"ottcher function is not defined in a unique way; however, if $\beta_A$ is some B\"ottcher function, then any other 
B\"ottcher function has the form $\beta_A\circ \v z,$ where $
\v^{n-1}=1.$ 
In this notation, our 
 main result concerning equation \eqref{abo} is the following.

\bt \l{mt1} Let $A\in \Gamma$ be a formal power series of order $n$, and 
$\beta_A$ some B\"otcher function. Then every 
 decomposition \be \l{eq1} A= A_1\circ A_2 \circ \dots \circ A_r\ee of $A$ into a composition 
of elements $A_1, A_2, \dots ,A_r$ of $\Gamma$ is equivalent to the decomposition 
\be \l{eq3} A=(\beta_A\circ  z^{\ord A_1})\circ z^{\ord A_2} \circ \dots \circ (z^{\ord A_r}\circ \beta_A^{-1}).\ee
Thus, equivalence classes 
 of decompositions of $A$  
are in a one-to-one correspondence with ordered factorizations of $n$. 
\et 

The main motivation for writing this paper was to construct in the formal power series setting  an analogue of the decomposition theory of  {\it rational functions}. Correspondingly, the definition of the equivalency of decompositions of 
elements of $\Gamma$ given above mimics the corresponding definition 
from the decomposition theory of rational functions, in which two decompositions  \eqref{decc} of a rational function of degree at least two $A$   into compositions of rational functions of degree at least two $A_1, A_2, \dots, A_k$ and $\widehat A_1, \widehat A_2, \dots, \widehat A_m$   are considered as equivalent if \eqref{decc0} holds for some  M\"obius transformations $\mu_i,$ $1\leq i \leq k-1.$  As expected, 
the results obtained in this paper differ significantly from the corresponding results for rational functions, generally being simpler. For instance, even for polynomial decompositions, the analogue of Theorem \ref{mt1}, obtained by Ritt (\cite{r1}), is substantially more complex. On the other hand, for arbitrary rational functions, such an analogue is not known, and typical results in the area primarily concern either decompositions of specific types of functions or functional equations of a particular form (see e.g.  \cite{az}, \cite{bog}, \cite{ere}, \cite{fz}, \cite{mp}, \cite{ngw}, \cite{lau}, \cite{semi}, \cite{rev}, \cite{ritt}).

The main method in the study of decompositions of rational functions is the monodromy method, which involves examining the monodromy group associated with a given rational function. On the other hand, the primary technical tool in the study of equations in formal power series is the Böttcher functions.   Our approach consists in the systematic use  along with the B\"ottcher functions what we call the {\it transition functions}. By definition, the transitions functions for $A\in \Gamma$ are elements $\phi_A$ of $k_1[[z]]$ satisfying  $$A\circ \phi_A=A.$$ For $A\in \Gamma$ of  order $n$ there exist exactly $n$ transition functions  forming a cyclic group with respect to the operation of composition. We will call this group the {\it transition group} and denote it by $G_A$. Although the transition groups are quite simple from a group-theoretic perspective, they turn out to be very convenient for studying equations in 
$\Gamma$  since the relative position of these groups within $k_1[[z]]$ reflects the mutual compositional properties of the corresponding elements of $\Gamma$. We illustrate this statement with the following two results, which we consider among the main results of the paper. 

The first result concerns the functional equation $F=X\circ A,$ where $F,A$ are given and $X$ is unknown. 

\bt \l{mt2}  Let $A\in k_n[[z]]$, $n\geq 2$, and $F\in k_{nm}[[z]]$, $m\geq 1.$  
 Then the equation 
$$ F=X\circ A$$ has a solution in $X\in k_m[[z]]$ if and only if $G_A\subseteq G_F.$ In particular, 
for $A,B\in \Gamma$  of the same order the equality $ G_A=G_B$ holds if and only if 
$  B=\mu\circ A$ for some $\mu\in k_1[[z]]$
\et

The second result concerns the functional equation $X\circ A=Y\circ B,$ where $A,B$ are given and $X,Y$ are unknown.

\bt \l{mt3} Let $A,B \in \Gamma$.  Then the equation  $$ X\circ A=Y\circ B$$ 
 has a solution in $X,Y\in zk[[z]]$ if and only if 
$$\phi_A\circ \phi_B=\phi_B\circ \phi_A $$ for all $\phi_A\in G_A$ and  $\phi_B\in G_B.$
\et 

Along with  decompositions of general elements of $\Gamma$, we study decompositions of elements of a special form.  Specifically, 
 we address  the following problem posed by Horwitz and Rubel in \cite{hr}: if $h$ is the composition of two formal power series $f$ and $g,$ and if $h$ is even, what can be said about $f$ and $g$?  
Some partial results on this problem and its modifications, concerning decompositions of entire functions or polynomials, were obtained in the papers \cite{bear0},  \cite{bear},  \cite{hr},  \cite{hor}. 

In this paper, we provide a complete solution to the problem of Horwitz and Rubel 
in the case where $h$ and $f,g$  are elements of $\Gamma.$ In fact, along with {\it even} formal power series, that is, series having the form $R(z^2)$ for some $R\in k[[z]]$, we also consider {\it odd} series having the form  $zR(z^2)$  and, more generally, {\it symmetric} series having the form 
$z^rR(z^m)$, where $m\geq 2$, $r\geq 0$ are integers.   Specifically, we prove the following result.

\bt \l{dee}  Let $A\in \Gamma$ be a formal power series of the form $A=z^rR(z^m)$, where $R\in k[[z]]$ and  $m\geq 2$,  $r\geq 0$  are integers. Then for any  decomposition $A=A_1\circ A_2$,
where $A_1,A_2\in \Gamma$, there exist $\mu\in k_1[[z]]$ and $R_1,R_2\in k[[z]]$ such that 
$$ A_1=z^{r_1}R_1(z^{\frac{m}{\gcd(r_2,m)}})\circ \mu^{-1}, \ \ \ \  A_2=\mu \circ z^{r_2}R_2(z^m)$$
for some integers $r_1, r_2\geq 0$  
satisfying the condition $r_1r_2\equiv r\, (\mod m).$
\et

Notice that Theorem \ref{dee} implies that if $A=A_1\circ A_2$ is even, then either $A_2$ is even, or there exists $\mu\in k_1[[z]]$ such that $\mu^{-1}\circ A_2 $ is odd and $A_1\circ \mu$ is even. On the other hand, 
if $A=A_1\circ A_2$ is odd, then Theorem \ref{dee}  implies that there exists $\mu\in k_1[[z]]$ such that $A_1\circ \mu$ and $\mu^{-1}\circ A_2$ are both odd (see Corollary \ref{bc2}).

As an application of our results about functional equations in $\Gamma$, we provide 
a handy necessary condition for a subsemigroup of $\Gamma$ to be {\it right amenable}, 
 meaning that it  
admits a finitely additive probability measure $\mu$ defined on all  subsets 
of $S$  
such that for all $a\in S$ and $T\subseteq S$ 
the equality $$ \mu(Ta^{-1}) = \mu(T)$$ holds, where   
 the set $Ta^{-1}$ is defined by the formula  $$ Ta^{-1}=\{s \in S\, | \,sa \in T\}.$$ 

Let us denote by $\f Z$ the  subsemigroup of $\Gamma$  consisting of  monomials  $az^n,$ where $a\in k^*$ and $n\geq 2,$ and by  $\f Z^U$ the subsemigroup consisting of all monomials  of 
the form $\omega z^n,$ $n\geq 2,$ where $\omega$ is a root of unity.
We say that two subsemigroups $S_1$ and $S_2$ of $\Gamma$ are {\it conjugate} if 
there exists a formal power series $\alpha\in k_1[[z]]$ such that $$\alpha\circ S_1\circ \alpha^{-1}=S_2.$$  
It was shown in \cite{amf} that a {\it finitely} generated subsemigroup of $\Gamma$ is right amenable if and only if it is conjugate to a subsemigroup of $\f Z^U$. 
However, it was observed that 
an {\it infinitely} generated  right amenable  subsemigroup of $\Gamma$ is not necessarily conjugate to a  subsemigroup of $\f Z^U$.
In this paper, we prove the following result. 

\bt \l{mt4} Every right amenable subsemigroup  $S$  of $\Gamma$ is conjugate to a subsemigroup of $\f Z.$ 
\et 

 Moreover, we show that the conclusion of Theorem \ref{mt4} holds already under the  assumption that $S$ is {\it right reversible}, which is a weaker condition than the assumption that $S$ is right amenable (see Theorem \ref{mt5}).  We deduce these results   from the following statement of independent interest.

\bt \l{51} Let $A,B \in \Gamma$  be formal power series, and 
$\beta_A$, $\beta_B$ some B\"otcher functions. Then  the equation 
\be \l{ggg} X\circ A^{\circ l}=Y\circ B^{\circ s}\ee has a solution in $X,Y\in zk[[z]]$
 for all $s,l\geq 1$ if and only if $\beta_A=\beta_B\circ cz$ for some $c\in k^*.$ 
\et 

Notice that Theorem \ref{51} includes the characterization of commuting elements of $\Gamma$ in terms of 
their B\"ottcher functions, as obtained by Dorfer and Woracek (\cite{ew}). Specifically, it implies that $A,B\in \Gamma$  commute if and only if $\beta_A=\beta_B\circ \v z$ for some $\v$ satisfying $$\v^{(\ord A-1)(\ord B-1)}=1$$ (see Corollary \ref{wor}). 

 This paper is organized as follows. In the second section, after recalling several elementary facts about the semigroup $k[[z]]$ we discuss B\"otcher functions and some of their immediate applications to functional equations. In the third section, we introduce transition functions and establish their basic properties. 
In the fourth section, we solve the functional equations  $$F=A\circ X \ \ \  {\rm and}  \ \ \ F=X\circ A,$$ where $F, A\in \Gamma$ are given  and $X\in zk[[z]]$ is unknown, in terms of the corresponding B\"ottcher functions. We also prove Theorem \ref{mt2} and several of its corollaries. 

In the fifth section, we apply the obtained results to decompositions of elements of $\Gamma$, and prove Theorem \ref{mt1}. In the sixths section, we characterize symmetric series  in terms of their B\"otcher and transition functions, and prove Theorem \ref{dee}. 
 We also reprove the result of Reznick (\cite{rez}) stating that if an iterate of $A\in \Gamma$  is symmetric, then $A$ is also symmetric. 
In the seventh section, we consider the functional equation $$X\circ A=Y\circ B,$$ where 
$A,B\in \Gamma$ are given and $X,Y\in zk[[z]]$ are unknown, and prove Theorem \ref{mt3} and 
Theorem \ref{51}. Finally, we establish the  aforementioned necessary condition for the right amenability and the right reversibility of subsemigroups of $\Gamma$.

\section{B\"ottcher functions}

\subsection{Lemmata about formal power series}
In this paper, $k$ always denotes an algebraically closed field of characteristic zero. 
Notice that  the number of $n$th roots of unity in such $k$ equals $n$ for every $n\geq 1$. We will denote by $U_n$
the group of $n$th roots of unity in $k$, and by $U_n^P$ the subset of $U_n$ consisting of primitive $n$th roots of unity. 

For elementary properties of the ring of formal power series $k[[z]]$ and the semigroup $zk[[z]]$ under the composition  operation  $\circ $, 
we refer the reader to the first paragraph of \cite{car}.  In particular, we will use the fact that  $k[[z]]$ is an integer domain and that 
an element $A$ of $zk[[z]]$ is invertible with respect to $\circ$  if and only if $A$ belongs to $k_1[[z]].$ 
Below we collect some further simple facts about $k[[z]]$. 

\bl \l{0} Formal power series $\mu_1,\mu_2\in k[[z]]$ satisfy the equality    
$$ z^n\circ \mu_1=z^n\circ \mu_2, \ \ \ n\geq 2,$$ if and only if 
$\mu_1=\v\mu_2$ for some $\v\in U_n.$
\el
\pr 
Since 
$$\mu_1^n-\mu_2^n=\prod_{\v \in U_n}
(\mu_1- \v\mu_2),$$ 
the lemma follows from the fact that $k[[z]]$ is an integer domain.  \qed 

\bl \l{00} Let  $\mu\in k[[z]]\setminus k$ and $a,b\in k^*$ 
  satisfy the equality    
\be \l{oo} \mu\circ  az=b z\circ \mu.\ee Then  $b=a^r$ for some $r\in \N$. Furthermore,  either  $\mu=cz^r$, $r\geq 1$, for some $c\in  k^*$,  or  $a$ is a root of unity. 
 Finally, 
 $\mu$ satisfies the equality    
\be \l{o} \mu\circ \v z={\v}^{\,r} z\circ \mu\ee 
 for some  $\v\in U_n^P$ and 
$r$, $0\leq r \leq n-1,$  if and only if there exists a
 formal power series $R\in  k[[z]]$ such that 
$\mu=z^rR(z^n).$
\el
\pr The proof is obtained by a comparison of coefficients in the left and the right parts of \eqref{oo} and \eqref{o}. \qed

\bl \l{11} A formal power series $\mu\in k[[z]]$  satisfies the equality    
\be \l{int} z^n\circ \mu=\mu\circ z^n, \ \ \ n\geq 2,\ee if and only if
$\mu=\v z^m$ for some $\v\in U_{n-1}$ and  $m\geq 0.$ 
\el
\pr Setting $m=\ord\, \mu$ and substituting  $\mu=\sum_{i=m}^{\infty}c_iz^i$ into
\eqref{int} we see that     $c_m^n=c_m.$ Furthermore, 
if  $\mu\neq c_mz^m$ we obtain a contradiction as follows. Let 
$l> m$ be the minimum number such that $c_l\neq 0$. Then  
$$\mu=c_{m}z^{m}+c_{l}z^{l}+{\rm higher \ \ terms},$$ implying that  
 $$\mu\circ z^n=c_{m}z^{mn}+c_lz^{ln}+{\rm higher \ \ terms}.$$
On the other hand, 
$$z^n\circ \mu=c_{m}^nz^{mn}+nc_{m}^{n-1}c_lz^{m(n-1)+l}+{\rm higher \ \ terms}. 
$$
Since $$m(n-1)+l<l(n-1)+l=ln,$$ this is impossible, and hence  $\mu= c_mz^m$.
 \qed

\bl \l{12} Formal power series $\mu_1,\mu_2\in k[[z]]\setminus k$ 
satisfy the equality    
\be \l{eq} z^n\circ \mu_1=\mu_2\circ z^n, \ \ \ n\geq 2,\ee if and only if there exist  $R\in  k[[z]]$ and $r$, $0\leq r\leq n-1,$   such that  
$$\mu_1=z^rR(z^n), \ \ \ \ \mu_2=z^rR^n(z).$$ 
\el
\pr The identity 
\be \l{ide} z^n\circ z^rR(z^n)=z^rR^n(z)\circ z^n\ee is checked 
 by a direct calculation. To prove the ``only if'' part, we observe that 
for any  $\v_n\in U_n^P $ 
equality $\eqref{eq}$ implies the equality  
$$z^n\circ \mu_1=z^n\circ (\mu_1 \circ \v_n z).$$ Therefore, by Lemma \ref{0},  there exists $r,$ $0\leq r \leq n-1,$ such that  $$\mu_1 \circ \v_n z=\v_n^r z \circ \mu_1,$$ implying by Lemma \ref{00} that 
$\mu_1=z^rR(z^n)$ for some $R\in k[[z]]$. 
 It follows now from \eqref{eq} that 
$$\mu_2\circ z^n=z^n\circ \mu_1=z^{rn}R^n(z^n)=z^{r}R^n(z)\circ z^n,$$ 
implying that $\mu_2=z^rR^n(z).$
\qed 

\vskip 0.2cm
Notice that the representation $\mu_2=z^rR^n(z)$ appearing in Lemma \ref{12} defines the series $R$ only up to a multiplication by an $n$th root of unity. Accordingly,  to $\mu_2$  
correspond $n$ different $\mu_1$ such that \eqref{eq} holds.  

\subsection{B\"ottcher functions and the equation $A\circ X=Y\circ B$}

Let $A\in \Gamma$ be 
 a formal power series of order $n$.  Then the corresponding B\"ottcher function is defined as a formal power series  $\beta_A\in k_1[[z]]$ 
such that  the equality 
\be \l{a} A\circ \beta_A=\beta_A\circ z^n\ee holds. 
 It is well known that such a function exists and  is defined in a unique way 
up  to the change $\beta_A\rightarrow \beta_A \circ \v z,$ where $\v\in U_{n-1}.$ 
In the context of complex dynamics, this fact is widely used and goes back to  B\"ottcher (see \cite{bo}, \cite{rr}, \cite{mil}). For the proof in the algebraic setting, see \cite{kau} (Hilffsatz 4). 
Notice that the map 
\be \l{map} A\rightarrow  \beta_X^{-1} \circ A \circ \beta_X, \ee where  $X$ is a fixed element of $\Gamma$,  is a semigroup automorphism of $\Gamma$.

Among other things, the existence of B\"ottcher functions yields the following statement.

\bt \l{buri} 
Let $A_1,A_2\in k[[z]]$ and  $X\in zk[[z]]$ be formal power series.  Then 
 the equality 
\be \l{b0} 
A_1\circ X=A_2\circ X
\ee
holds  if and only $A_1=A_2.$ 
\et 
\pr In case $X$ is invertible in the semigroup $zk[[z]]$, the statement is clear. Otherwise setting $n=\ord X$ and conjugating  \eqref{b0} by $\beta_X$, we obtain the equality 
$$(\beta_X^{-1} \circ A_1 \circ \beta_X)\circ z^n=(\beta_X^{-1} \circ A_2 \circ \beta_X)\circ z^n, 
$$
which obviously implies that 
$$\beta_X^{-1} \circ A_1 \circ \beta_X=\beta_X^{-1} \circ A_2 \circ \beta_X.$$ Since \eqref{map} is an isomorphism, this implies in turn that $A_1=A_2$. 
\qed

\vskip 0.2cm 

Using B\"ottcher functions, one can provide a solution in  $X,Y\in zk[[z]]$ of the functional equation $$A\circ X=Y\circ B,$$ where $A$ and $B$ are given elements of $\Gamma$ of the same order,     generalizing equation \eqref{a}. 
We start by considering the following particular case. 
\bt \l{21}
Let $A\in \Gamma$ be a formal power series of order $n$, and  $\beta_A$ some B\"ottcher function. 
 Then solutions of the equation 
\be \l{b} 
A\circ X=Y\circ z^n
\ee
 in $X,Y\in zk[[z]]$ are given by the formulas 
\be \l{kor} X=\beta_A\circ z^r R(z^n), \ \ \ \ Y=\beta_A\circ z^r R^n(z),\ee 
where $R\in  k[[z]]$ and $0\leq r\leq n-1.$
Furthermore, if $X=Y$, then solutions of \eqref{b}  are given by the formula 
\be \l{lol} X=\beta_A\circ \v z^l,\ \ \ \ \v \in U_{n-1},\ee where $l=\ord X$. 
\et 
\pr 
The fact that $X$ and $Y$ defined by \eqref{kor} satisfy \eqref{b} follows from equalities \eqref{ide} and \eqref{a}. On the other hand, if  \eqref{b} holds, then 
 taking an arbitrary B\"ottcher function $\beta_A$ and substituting   
$\beta_A\circ z^n\circ \beta_A^{-1}$ for $A$ in \eqref{b}, we obtain 
$$\beta_A\circ z^n\circ \beta_A^{-1}\circ X=Y\circ z^n,$$
implying that 
$$ z^n\circ (\beta_A^{-1}\circ X)=(\beta_A^{-1}\circ Y)\circ z^n.$$ Thus,  equalities \eqref{kor} hold  by Lemma \ref{12}. 

Furthermore, if $X=Y,$ then \eqref{kor} implies that 
$$z^rR(z^n)=z^rR^n(z).$$ In turn, this yields that $R$ commutes with 
$z^n$, implying by Lemma \ref{11} that $R=\v z^{m}$, where $\v\in U_{n-1}$ and $m\geq 0.$ Therefore, 
$$X=z^rR(z^n)=\v z^{l},$$ where $$l=\ord z^rR(z^n)=\ord X. \eqno{\Box}$$
\vskip 0.2cm

Theorem \ref{21} implies the following more general statement. 
\bt \l{211}
Let $A,B\in \Gamma$ be formal power series of the same order  $n$, and $\beta_A$, $\beta_B$ some B\"ottcher functions. 
Then solutions of the equation 
\be \l{bb} 
A\circ X=Y\circ B
\ee
in $X,Y\in zk[[z]]$ are given by the formulas  
\be \l{korb} X=\beta_A\circ z^r R(z^n) \circ \beta_B^{-1} , \ \ \ \ Y=\beta_A\circ z^r R^n(z) \circ \beta_B^{-1}, \ee  
where $R\in  k[[z]]$ and $0\leq r\leq n-1.$ 
Furthermore, if $X=Y$, then solutions of \eqref{bb}  are given by the formula 
\be \l{lol1} X=\beta_A\circ \v z^l \circ \beta_B^{-1},\ \ \ \ \v\in U_{n-1},\ee where    $l=\ord X$. 
\et 
\pr  For an arbitrary  B\"ottcher function $\beta_B$,  
equality \eqref{bb} is equivalent to the equality 
$$A\circ (X\circ \beta_B)=(Y\circ \beta_B)\circ z^n.$$  
Thus, the theorem follows from Theorem \ref{21}. \qed

\section{Transition  functions}

Let $A\in \Gamma$ be 
 a formal power series of order $n$.  We recall that we defined  
 transition functions for $A$  as formal series $\phi_A$ 
satisfying 
\be \l{c} A\circ \phi_A=A. \ee
 It is clear that such series necessarily belong to $k_1[[z]]$ and form a group, which we denote by $G_A$.

The following two lemmas are modifications of the results of Section 2 in \cite{han} characterizing solutions of \eqref{c} in the analytical setting.

\bl \l{31} Let 
 $A\in \Gamma$ be a formal power series, and $\beta_A$ some  B\"ottcher function. Then 
\be \l{fina} G_A=\{\beta_A \circ \v z\circ \beta_A^{-1}\ \vert\   \v\in U_n\}.\ee

\el 
\pr It follows from equality \eqref{a} that for every $\v\in U_n$ we have  
$$A\circ \beta_A=A\circ \beta_A \circ \v z,$$  implying that 
$$A=A\circ (\beta_A \circ \v z\circ \beta_A^{-1}).$$

On the other hand,  if equality  \eqref{c} holds, then conjugating its parts by $\beta_A$, we obtain  
$$z^n\circ (\beta_A^{-1} \circ  \phi_A\circ \beta_A)=z^n,$$ implying by Lemma \ref{0} that 
$\beta_A^{-1} \circ  \phi_A\circ \beta_A=
\v z$ for some $\v\in U_n$. \qed

\vskip 0.2cm

For a formal power series $\phi\in k_1[[z]]$, we denote by $|\phi|$  the order of $\phi$ in the group $k_1[[z]]$. Thus, $|\phi|$ equals the minimum number $d$ such that 
$\phi^{\circ d}=z,$ if such a number exists, and  $|\phi|$ equals $\infty$, if $\phi^{\circ d}$ is distinct from $z$ for every $d\geq 1$.

\bl \l{34}  
Let $\phi\in k_1[[z]]$ 
 be a formal power series with $\vert \phi\vert=d$. Then $\phi=\phi_A$ for some formal power series $A\in \Gamma$  if and only if 
 $1<d<\infty$. Moreover, in the last case $\phi=\phi_A$ for some  $A$ of order $d$. 
\el 
\pr Since the functions defined by  \eqref{fina} satisfy  $\phi_A^{\circ n}=z,$ the ``only if'' part follows from Lemma \ref{31}. On the other hand, if $1<d<\infty,$ then 
 setting 
\be \l{aaa} A=z\cdot \phi\cdot \phi^{\circ 2}\cdot\, \dots \,\cdot \phi^{\circ (d-1)},\ee we see that 
$A\in k_d[[z]]$ and 
the equality $A\circ \phi =A$ holds. \qed 

\vskip 0.2cm

The following lemma  follows immediately from Lemma \ref{31}. 

\bl \l{if} Let $A\in \Gamma$. Then $G_A$ is a cyclic group of order $n$, whose generators are $\beta_A \circ \v_n z\circ \beta_A^{-1},$ where $\v_n\in U_n^P$. \qed 
\el

The following lemma relates the transition group  for $A\in \Gamma$ with the transition groups for 
$A^{\circ l},$ $l\geq 1,$ and 
$$A_{\mu}=\mu^{-1} \circ A \circ \mu,\ \ \  \mu\in k_1[[z]].$$

\bl \l{fif} Let $A\in \Gamma$ be a  formal power series of order $n$, and $\beta_A$ some  B\"ottcher function. Then  
\be \l{go1}  G_{A^{\circ l}}=\{\beta_A \circ \v z\circ \beta_A^{-1}\ \vert\   \v\in U_{nl}\}, \ \ \ \ l\geq 1,\ee and  
\be \l{conj} G_{A_{\mu}}=\mu^{-1} \circ G_A \circ \mu, \ \ \ \ \mu\in k_1[[z]].\ee
\el
\pr Equality \eqref{go1} follows from Lemma \ref{31} and the fact that $\beta_A$ remains a B\"ottcher function for $A^{\circ l}$, $l\geq 1.$  On the other hand, since $\ord A=\ord A_{\mu}=n,$ 
equality \eqref{conj} follows from the equality $$ A_{\mu}\circ (\mu^{-1} \circ \phi_A \circ \mu)=A_{\mu},\ \ \ \phi_A\in G_A,$$ 
 which is obtained by a direct calculation. \qed

The following statement is a counterpart of Theorem \ref{buri} for the functional equation  
$A\circ X_1=A\circ X_2$.

\bt \l{41} 
Let $A\in \Gamma$  and $X_1,X_2\in zk[[z]]$. Then  the equality 
\be \l{d} A\circ X_1=A\circ X_2\ee 
holds if and only if 
$$X_2=\phi_A\circ X_1$$ for some $\phi_A\in G_A.$
\et 
\pr The ``if'' part is obvious. On the other hand if equality \eqref{d} holds, then 
conjugating its parts by $\beta_A$ 
we obtain 
$$z^n\circ (\beta_A^{-1} \circ  X_1\circ \beta_A) =z^n\circ  (\beta_A^{-1} \circ  X_2\circ \beta_A),$$ 
implying that  
 $$\beta_A^{-1} \circ  X_2\circ \beta_A=\v z\circ \beta_A^{-1} \circ  X_1\circ \beta_A$$ 
for some $\v \in U_n$ by Lemma \ref{0}.  
Therefore, 
  $$ X_2 =\beta_A \circ\v z\circ \beta_A^{-1}\circ X_1=\phi_A\circ X_1 $$
 by Lemma \ref{31}.
\qed 

\vskip 0.2cm

\section{Functional equations $F=A\circ X$ and $F=X\circ A$} 
The next two results provide solutions of the functional equations  $F=A\circ X$ and $F=X\circ A$, where $F,A\in \Gamma$ are given  and $X\in zk[[z]]$ is unknown, in terms of the corresponding B\"ottcher functions 
$\beta_F$ and $\beta_A.$ 

\bt \l{231} Let $A\in k_n[[z]]$, $n\geq 2$, and $F\in k_{nm}[[z]]$, $m\geq 1,$ be formal power series,  and 
$\beta_A$,  $\beta_F$ some 
B\"ottcher functions. Then the equation 
\be \l{eqa} F=X\circ A\ee
has a solution in $X\in k_m[[z]]$ if and only if there exist $R\in k[[z]]$ and $r,$ $0\leq r \leq n-1,$  such that 
\be \l{eqa1} z^m\circ \beta^{-1}_F\circ \beta_A=z^rR(z^n).\ee
 Furthermore, if \eqref{eqa1} holds, then \eqref{eqa} has a unique solution $X$ given by the formula  
 \be \l{eqa2} X=\beta_F\circ z^rR^n(z)\circ \beta_A^{-1}.\ee

\et 
\pr 
Substituting  
$\beta_F\circ z^{nm}\circ \beta_F^{-1}$ for $F$ and   
$\beta_A\circ z^n\circ \beta_A^{-1}$ for $A$ to  \eqref{eqa}, we obtain 
the equality 
$$\beta_F\circ z^{nm} \circ \beta_F^{-1}=X\circ \beta_A\circ z^{n} \circ \beta_A^{-1},$$ which in turn implies the equality 
$$z^{n} \circ (z^{m} \circ \beta_F^{-1}\circ \beta_A)=(\beta_F^{-1}\circ X\circ \beta_A)\circ z^{n}.$$
Hence, the ``only if'' part follows from Lemma \ref{12}. 

In the other direction,  \eqref{eqa1} implies that 
$$F=\beta_F\circ z^{nm} \circ \beta_F^{-1}=\beta_F\circ z^{n} \circ z^{m}\circ \beta_F^{-1}=\beta_F\circ z^{n} \circ  z^rR(z^n)\circ \beta_A^{-1}=$$
$$=\beta_F\circ z^rR^n(z)\circ z^n\circ \beta_A^{-1}=\beta_F\circ z^rR^n(z)\circ \beta_A^{-1}\circ A.
$$
Thus, \eqref{eqa} holds for $X$ given by \eqref{eqa2}.  
Finally, the function $X$ is defined by formula \eqref{eqa2} in a unique way by Theorem \ref{buri}.
\qed

\bt \l{2mt2} Let $A\in k_n[[z]]$, $n\geq 2$, and $F\in k_{nm}[[z]]$, $m\geq 1,$ be formal power series, and 
$\beta_A$,  $\beta_F$ some 
B\"ottcher functions.  Then the equation 
\be \l{2eqa} F=A\circ X\ee
has a solution in $X\in k_m[[z]]$ if and only if 
there exist $L\in k[[z]]$ and $r,$ $0\leq r \leq n-1,$  such that 
\be \l{2eqa1} \beta^{-1}_A\circ \beta_F \circ z^m =z^rL^n(z).\ee 
 Furthermore, if  \eqref{2eqa1} holds, then \eqref{2eqa} has $n$ solutions given by the formula 
 \be \l{2eqa2} X=\beta_A\circ \v z \circ z^r L(z^n)\circ \beta_F^{-1}, \ \ \ \ \v\in U_n.\ee

\et 
\pr Equality \eqref{2eqa} implies the equality 
$$\beta_F\circ z^{nm} \circ \beta_F^{-1}=\beta_A\circ z^{n} \circ \beta_A^{-1}\circ X,$$ which in turn implies the equality 
$$(\beta_A^{-1}\circ \beta_F\circ z^{m})\circ z^{n} =z^{n}\circ (\beta_A^{-1}\circ X\circ  \beta_F).$$
Therefore, the ``only if'' part follows from Lemma \ref{12}. 

In the other direction,  \eqref{2eqa1} implies that 
$$F=\beta_F\circ z^{nm} \circ \beta_F^{-1}=\beta_F\circ z^{m} \circ z^{n}\circ \beta_F^{-1}=\beta_A \circ z^rL^n(z)\circ z^{n} \circ  \beta_F^{-1}=$$
$$=
\beta_A\circ z^n\circ z^rL(z^n)\circ  \beta_F^{-1}=A \circ \beta_A\circ z^rL(z^n)\circ \beta_F^{-1}.
$$
Thus, \eqref{2eqa} holds for $$X=\beta_A\circ z^rL(z^n)\circ \beta_F^{-1}.$$ 
Finally, by Theorem \ref{41} and Lemma \ref{31}, any other solution of \eqref{eqa2} has the form 
$$X=\phi_A\circ \beta_A\circ z^rL(z^n)\circ \beta_F^{-1}=\beta_A \circ \v z\circ \beta_A^{-1}\circ \beta_A\circ z^rL(z^n)\circ \beta_F^{-1}=$$
$$=\beta_A\circ \v z\circ z^rL(z^n)\circ \beta_F^{-1}, \ \ \ \ \v\in U_n. \eqno{\Box}$$

\vskip 0.2cm

\noindent {\it Proof of Theorem \ref{mt2}}. 
 If 
$ F=X\circ A,$ then for any $\phi_A\in G_A$ we have 
$$F\circ \phi_{A}=X\circ A\circ \phi_{A}=X\circ A=F,$$
implying that $G_A\subseteq G_{F}.$

In the other direction, the equality $F\circ \widehat\phi_{A}=F$ for some  generator 
 $\widehat\phi_A$ of $G_A$ implies that  
\be \l{it} \beta_F\circ z^{nm} \circ \beta_F^{-1}\circ \beta_A\circ  \v_{n}z \circ \beta_A^{-1}=\beta_F\circ z^{nm} \circ \beta_F^{-1}\ee  
  for some B\"ottcher functions $\beta_A$, $\beta_F$ and 	
	 $\v_n\in U_n^P$.  It is clear that equality \eqref{it} implies the equalities 
$$z^{nm} \circ \beta_F^{-1}\circ \beta_A\circ  \v_{n}z=z^{nm} \circ \beta_F^{-1}\circ \beta_A$$ and 
$$z^n\circ (z^{m} \circ \beta_F^{-1}\circ \beta_A\circ  \v_{n}z)=z^n\circ (z^{m} \circ \beta_F^{-1}\circ \beta_A).$$ 
In turn, the last equality implies by Lemma \ref{0} that 
$$(z^{m} \circ \beta_F^{-1}\circ \beta_A)\circ  \v_{n}z=\v_n^rz\circ  (z^{m} \circ \beta_F^{-1}\circ \beta_A)$$  for some   $r,$ $0\leq r \leq n-1.$ 
It follows now from Lemma \ref{00} that 
there exists   $R\in k[[z]]$ such that  \eqref{eqa1} holds. Therefore, the equality $ F=X\circ A$ holds for   some  $X\in k_m[[z]]$ 
by Theorem \ref{231}.
\qed

\vskip 0.2cm
For brevity, we will say that $A\in \Gamma$ is a {\it compositional right factor} of $F\in \Gamma$ if there exists $X\in zk[[z]]$ such that  
 $F=X\circ A$. Compositional left factors are defined similarly.

\bc \l{iit} 
 Let $F\in \Gamma$ be a formal power series,   and $A,B\in \Gamma$ some  compositional right factors of $F$. Then any $\phi_A\in G_A$ and $\phi_B\in G_B$ commute.   
\ec
\pr By Theorem \ref{mt2}, any $\phi_A\in G_A$ and $\phi_B\in G_B$  are elements of the commutative group $G_F.$  \qed

\vskip 0.2cm
The following corollary provides a criterion for two elements of $\Gamma$ to have  a ``common'' compositional right factor in $\Gamma.$

\bc \l{35}  Let $A \in k_n[[z]]$, $B\in k_m[[z]]$, $n,m\geq 2$,  be formal power series, and $d\geq 2$  a common divisor of $n$ and $m$.  
Then the system 
\be \l{te} 
A= \widetilde  A\circ W, \ \ \ \ B=  \widetilde  B\circ W, 
\ee 
has a solution in $ \widetilde   A\in k_{n/d}[[z]]$, $ \widetilde  B\in k_{m/d}[[z]]$, and $W\in k_d[[z]]$
if and only if 
the intersection of the groups $G_A$ and $ G_B$ contains a group of order $d$.

\ec 
\pr Assume that  \eqref{te} holds and let $ \widehat\phi_W$ be a generator of $G_W$. Then by the ``only if'' part of  Theo\-rem \ref{mt2}  
\be \l{poi} \widehat\phi_W=\widehat\phi_A^{\circ n/d}=\widehat\phi_B^{\circ m/d} 
\ee
for some generator $\widehat\phi_A$ of $G_A$ and  some generator $\widehat\phi_B$ of $G_B$. 
Thus, $G_A\cap G_B$ contains a cyclic group of order $d$ generated by $\widehat\phi_W.$

In the other direction, if $G_A\cap G_B$ contains a  group of order $d$, and $\phi$ is its generator, then 
$$\phi=\widehat\phi_A^{\circ n/d}=\widehat\phi_B^{\circ m/d} 
$$
for some generator $\widehat\phi_A$ of $G_A$ and  some generator $\widehat\phi_B$ of $G_B$. On the other hand, 
since $\vert \phi \vert =d$, it follows from Lemma \ref{34} that  $\phi=\widehat\phi_W $ for some  $W\in k_d[[z]]$.  Using now the ``if'' part of Theorem \ref{mt2}, 
we conclude that \eqref{te} holds. \qed

\vskip 0.2cm
We finish this section by the following result,  providing  a criterion for a formal power series $D\in \Gamma$ to be a compositional right factor of a composition of formal power series $A,C\in \Gamma$.

\bt Let $A,C,D\in \Gamma$ be formal power series. Then the equation 
\be \l{fi} A\circ C=X\circ D \ee 
has a solution in $X\in k[[z]]$ if and only if for any $\phi_D\in G_D$ there exists $\phi_A\in G_A$ such that  
\be \l{gg} C\circ \phi_D=\phi_A \circ C. \ee
\et 
\pr If for any $\phi_D\in G_D$  equality \eqref{gg} holds for some  $\phi_A\in G_A$, then for any $\phi_D\in G_D$ we have 
$$A\circ C\circ \phi_D=A\circ \phi_A \circ C=A\circ C.$$
Therefore, $G_D\subseteq G_{A\circ C}$ and hence \eqref{fi} has a solution by Theorem \ref{mt2}. 

In the other direction, equality \eqref{fi} implies that 
$$A\circ C=A\circ C\circ \phi_D.$$
Thus, \eqref{gg} holds by Theorem \ref{41}. \qed

\section{Equivalency classes of decompositions of formal power series}

In this section, we prove Theorem \ref{mt1} and  deduce from it a corollary, which can be considered as  an analogue of the result of Engstrom (\cite{en}) about polynomial solutions of the equation $A\circ C=B\circ D.$

\vskip 0.2cm
\noindent{\it Proof of Theorem \ref{mt1}.}
Let 
\be \l{eq11} A=A_1\circ A_2 \circ \dots \circ A_r\ee 
be a decomposition of $A\in \Gamma$ with   $$\ord A_k=n_k, \ \ \ 1\leq k \leq r.$$ Since 
$$ \beta_A^{-1}\circ  A\circ \beta_A=z^n=(\beta_A^{-1}\circ  A_1)\circ A_2 \circ \dots \circ (A_r\circ \beta_A),$$ to prove the theorem it is enough to show that for $A=z^n$
   every  decomposition \eqref{eq11}  is equivalent to the decomposition 
\be \l{eq2} z^n= z^{n_1}\circ z^{n_2} \circ \dots \circ z^{n_r}.\ee
We prove the last statement by induction on $r$.  

Clearly,  $G_{z^n}=\{\v z \ \vert\   \v\in U_n\}$.  
Since $\vert G_{A_r}\vert =n_r$ and   $G_{A_r}$ is a subgroup of  $G_{z^n}$ by Theorem \ref{mt2}, this implies that  $G_{A_r}= \{\v z \ \vert\   \v\in U_{n_r}\}$. 
Thus, $G_{A_r}=G_{z^{n_r}}$, 
implying by Theorem \ref{mt2} that 
\be \l{z1} A_r=\mu_{r-1} \circ z^{n_r}\ee  for some $\mu_{r-1}\in k_1[[z]].$ 
Hence, if $r=2$,  we have  
$$z^{n_1n_2}=A_1\circ \mu_1 \circ z^{n_2},$$ implying by  Theorem \ref{buri} that 
$A_1=z^{n_1}\circ \mu_1^{-1}.$ On the other hand, if $r>2$, then in a similar way we obtain the equalities \eqref{z1} and \be \l{last} z^{n_1n_2\dots n_{r-1}}=A_1\circ A_2 \dots (A_{r-1}\circ \mu_{r-1}).\ee By the induction assumption,  
the decomposition in the right part of \eqref{last} is equivalent to the decomposition $z^{n_1}\circ z^{n_2} \circ \dots \circ z^{n_{r-1}}$, and in virtue of \eqref{z1} this implies that for $A=z^n$ 
every decomposition  \eqref{eq11}  is equivalent to decomposition \eqref{eq2}. \qed 
\vskip 0.2cm

\bc \la{eng}
Assume that  $A,B,C,D\in \Gamma$ satisfy \be \l{abc} A\circ C=B\circ D.\ee Then there exist 
$U, V, \widetilde A, \widetilde C, \widetilde B, \widetilde D\in zk[[z]], $  where
$$\ord U=\GCD(\ord A,\ord B),  \ \ \ \ord V=\GCD(\ord C,\ord D),$$
such that
\be \l{end} A=U\circ \widetilde A, \ \  B=U\circ \widetilde B, \ \ C=\widetilde C\circ V, \ \  D=\widetilde D\circ V,\ee
and \be \l{yba} \widetilde A\circ \widetilde C=\widetilde B\circ \widetilde D.\ee 
\ec
\pr Let us set $$F=A\circ C=B\circ D,$$ $$n=\ord F, \ \ a=\ord A, \ \ b=\ord B, \ \ c=\ord C, \ \ d=\ord D,$$  
$$u=\gcd(a,b), \ \ \ \  v=\gcd(c,d).$$   
Taking a B\"otcher  function $\beta_F$ and applying Theorem \ref{mt1}, we see that there exist $\nu,\mu\in k_1[[z]]$ such that 
$$A=\beta_F\circ z^a\circ \nu^{-1}, \ \ \ C=\nu \circ z^c\circ \beta_F^{-1},$$ and 
$$B=\beta_F\circ z^b\circ \mu^{-1}, \ \ \ D=\mu \circ z^d\circ \beta_F^{-1}.$$
Therefore, the statement of the corollary is true for 
$$U=\beta_F\circ z^u, \ \ \ \  V=z^v\circ \beta_F^{-1}$$ and 
$$
\widetilde A= z^{\circ\frac{a}{u}}\circ \nu^{-1}, \ \ \widetilde C=\nu\circ z^{\circ\frac{c}{v}},\ \ \widetilde B= z^{\circ\frac{b}{u}}\circ \mu^{-1}, \ \ \widetilde D=\mu\circ z^{\circ\frac{d}{v}}.
$$

\section{Formal power series 
with symmetries}
\subsection{Characterizations of formal powers series 
with symmetries}

The following result characterizes 
elements of $\Gamma$ of the form $A=z^rR(z^m)$, where $R\in k[[z]]$ and $m\geq 2$, $r\geq 0$ are integers,   in terms of the corresponding B\"otcher  functions.

\bt \l{boee}  Let $A\in \Gamma$. Then $A$ has the form $A=z^rR(z^m)$ for some $R\in k[[z]]$ and 
integers  $m\geq 2$, $r\geq 0$  if and only if any B\"otcher  function $\beta_A$  has the form  $\beta_A=z L(z^m)$ for some $L\in k_0[[z]]$.
\et 
\pr Assume that for some B\"otcher function  $\beta_A$ the equality $\beta_A=z L(z^m)$ holds. Then
$\beta_A$ commutes with $\v_mz$ for any 
 $\v_m\in U_m^P$,   whence 
$$ (A\circ  \v_mz) \circ \beta_A=A\circ \beta_A\circ \v_mz=\beta_A\circ z^n\circ \v_mz=\beta_A\circ \v_{m}^nz \circ  z^n=$$ $$=\v_m^{n} z \circ\beta_A\circ  z^n=(\v_m^{n} z \circ A) \circ \beta_A.$$ Therefore, 
$$A\circ  \v_mz=\v_m^{n}z  \circ A,$$ implying by Lemma \ref{00} that  $A=z^rR(z^m)$.

In the other direction, let us assume that $A=z^rR(z^m)$ and 
set $\widehat A=z^rR^m(z).$ 
Since 
$$\widehat A\circ z^m=z^m\circ A,$$ for any B\"otcher  function $\beta_A$ we have 
$$\widehat A\circ (z^m\circ \beta_A)=z^m\circ A\circ \beta_A=(z^m\circ \beta_A)\circ z^n,$$ where $n=\ord A$, 
implying by Theorem \ref{21} that  
$$
z^m\circ \beta_A=\beta_{\widehat A}\circ \v z^m=(\beta_{\widehat A}\circ \v z) \circ z^m$$ for some  B\"otcher  function $\widehat \beta_A$  and $\v \in U_{n-1}$. 
 By Lemma \ref{12}, this implies that $\beta_A=z^l L(z^m)$, where $L\in k[[z]]$ and  $0\leq l\leq m-1$.  
Finally, 
since $\beta_A\in k_1[[z]],$ we conclude that 
$l=1$ and $L\in k_0[[z]].$ \qed 

\vskip 0.2cm

Notice  that if some B\"otcher  function  has the form  $\beta_A=z L(z^m)$, then all B\"otcher  functions have such a form. 

The following result is a counterpart of Theorem \ref{boee} in the context of transition functions.

\bt \l{traa}  Let $A\in \Gamma$. Then $A$ has the form $A=\mu \circ z^rR(z^m)$ for some  $\mu\in k_1[[z]]$, $R\in k[[z]]$, and integers  $m\geq 2$, $r\geq 0$  if and only if any transition function $\phi_A$ has the form $\phi_A=z M(z^m)$ for some $M\in k_0[[z]]$.
\et 
\pr Let us fix $\v_m\in U_m^P$. If some $\phi_A\in G_A$ has the form $\phi_A=z M(z^m)$, then $\phi_A$ commutes with $\v_mz$, implying that 
$$A\circ \v_mz=A\circ   \phi_A \circ \v_mz=(A\circ \v_mz)\circ   \phi_A.$$ Thus,  $ \phi_A$ belongs to $G_{A\circ \v_mz}.$ Therefore, if  any $\phi_A\in G_A$ has the above form, then $G_A=G_{A\circ \v_mz},$   implying by Theorem \ref{mt2} that 
\be \l{ew} A\circ \v_mz=\nu \circ A\ee for some $\nu\in k_1[[z]].$ 

Since \eqref{ew}  implies that $$A\circ (\v_mz)^{\circ l}=\nu^{\circ l} \circ A,\ \ \ l\geq 1,$$ the 
number $d=\vert \nu \vert$ 
is finite and divides $m$. If $d=1$, that is, if $\nu=z$, then applying Lemma \ref{00} to equality \eqref{ew} we conclude that $A=R(z^m)$ for some  $R\in k[[z]].$ On the other hand, if $d>1$, then 
 $\nu=\phi_{F}$ for some $F\in k_d[[z]]$
by Lemma \ref{34}, and hence 
$$\nu=\beta_F\circ \v z\circ \beta_F^{-1}$$ for some B\"otcher function  $\beta_F$ and $\v \in U_d$ by  
Lemma \ref{31}.  Moreover, since $d$ divides $m$, the equalities  $\v =\v_m^r$ and 
$$\nu=\beta_F\circ \v_m^rz\circ \beta_F^{-1}$$ hold for some $r,$ $0\leq r\leq m-1$.  
Substituting the right part of the last equality for $\nu$ in \eqref{ew}, we see that 
 $$(\beta_F^{-1} \circ A)\circ \v_mz= \v_m^rz\circ (\beta_F^{-1} \circ A).$$   
Hence,   by Lemma \ref{00}, 
$$\beta_F^{-1} \circ A=z^rR(z^m),$$ for some $R\in k[[z]]$. Thus,
the equality $A=\mu \circ z^rR(z^m)$ holds for $\mu=\beta_F$.

In the other direction, if $A=\mu \circ z^rR(z^m)$, then applying Corollary \ref{iit} to the function $$F=\widehat A\circ z^m=z^m\circ \mu^{-1} \circ A,$$ where  
 $\widehat A=z^rR^m(z)$,   
 we conclude that any $\phi_A\in G_A$ commutes with  $\phi_{z^m}=\v_mz$. Therefore, any $\phi_A$ has the form  $\phi_A=z M(z^m)$  by Lemma \ref{00}. \qed 
\vskip 0.2cm

\bc \l{cd} Let $A\in \Gamma$. Then $A$ has a compositional right factor $C\in \Gamma$ of the form $C= z^rR(z^m)$ for some   $R\in k[[z]]$ and integers  $m\geq 2$, $r\geq 0$  if and only if some transition function $\phi_A\neq z$ has the form $\phi_A=z M(z^m)$ for some $M\in k_0[[z]]$.
\ec 
\pr If $A$ has such a factor, then  by Theorem \ref{mt2} the group $G_A$ contains the non-trivial group 
$G_{C}$ as a subgroup. Moreover, all elements of the last group have the form $z M(z^m)$ by Theorem \ref{traa}.

In the other direction, let us assume that some transition function $\phi_A\neq z$ has the form $\phi_A=z M(z^m)$ and set  $d=\vert \phi_A\vert.$ By Lemma \ref{34}, $\phi_A=\phi_C$ for some $C\in \Gamma$ of order $d,$ and it is clear that $G_C=\langle \phi_A\rangle$. Thus, $A=B\circ C$ for some $B\in zk[[z]]$ by Theorem \ref{mt2}.  Moreover, since any iterate of a series of the form $z M(z^m)$ also has such form, it follows from $G_C=\langle \phi_A\rangle$ by Theorem \ref{traa} that $C$ has the form $\mu\circ z^rR(z^m)$ for some $\mu\in k_1[[z]]$. Finally, changing 
$B$ to $B\circ \mu$, we may assume that $C=z^rR(z^m)$. \qed

\subsection{Decompositions of formal powers series 
with symmetries}
Below, we provide some applications of Theorem \ref{boee} and Theorem \ref{traa}. 
    We start by proving Theorem \ref{dee}.

\vskip 0.2cm

\noindent {\it Proof of Theorem \ref{dee}}. Let us fix $\v_m\in U_m^P$.  Let 
\be \l{a1} A=A_1\circ A_2,\ee  be a decomposition of $A$ with $A_1,A_2\in \Gamma$.   
Considering the equality $$\widehat A\circ z^m=(z^m\circ A_1)\circ A_2,$$ where  
 $\widehat A=z^rR^m(z)$,   and using Corollary \ref{iit}, we see that any  $\phi_{A_2}\in G_{A_2}$ commutes with 
the  transition function $\phi_{z^m}=\v_mz$.   Thus,  any  $\phi_{A_2}\in G_{A_2}$ has the form $z M(z^m)$ for some $M\in k_0[[z]]$ by Lemma \ref{00}, and hence  
\be \l{a2} A_2=\mu \circ z^{r_2}R_2(z^m)\ee for some $\mu\in k_1[[z]]$, $R_2\in k[[z]]$,  and $r_2\geq 0$, by Theorem \ref{traa}.  

Furthermore, it follows from the equality 
$$z^rR(z^m)=A_1\circ \mu\circ z^{r_2}R_2(z^m)$$ that 
$$\big(A_1\circ \mu\circ z^{r_2}R_2(z^m)\big)\circ \v_m z=\v_m^{r}z\circ \big(A_1\circ \mu\circ z^{r_2}R_2(z^m)\big),$$ 
implying that 
$$A_1\circ \mu\circ \v_m^{r_2}z\circ z^{r_2}R_2(z^m)=\v_m^{r}z\circ A_1\circ \mu\circ z^{r_2}R_2(z^m)$$ and 
$$A_1\circ \mu\circ \v_m^{r_2}z=\v_m^{r}z\circ A_1\circ \mu.$$ 
 Since $\v_m^{r_2}$ is a primitive $\frac{m}{\gcd(r_2,m)}$th 
root of unity, it follows now from 
 Lemma \ref{00} that  
$$A_1\circ \mu=z^{r_1}R_1(z^{\frac{m}{\gcd(r_2,m)}})$$ 
 for some  $R_1\in k[[z]]$ and $r_1\geq 0$. 
 Thus, 
\be \l{a3} A_1=z^{r_1}R_1(z^{\frac{m}{\gcd(r_2,m)}})\circ \mu^{-1}. \ee 
Finally, it follows from \eqref{a1} and \eqref{a2}, \eqref{a3}  that $r_1r_2\equiv r\ (\mod m)$. \qed 

\vskip 0.2cm

Notice that in general  
the series $A_1$ in a decomposition $A=A_1\circ A_2$ of a symmetric series $A$ is ``less symmetric'' than $A$. Moreover,  
if $r_2=0$, then $A_1$ may be not symmetric at all. Nevertheless, the following statement is true.

\bc \l{bc1}  Let $A\in \Gamma$ be a formal power series of the form $A=z^rR(z^m)$, where $R\in k[[z]]$ and  $m\geq 2$,  $r\geq 1$  are  integers such that   $\gcd(r,m)=1.$ Then 
for any decomposition $A=A_1\circ A_2$, where   $A_1,A_2\in \Gamma$, there exist  $R_1,R_2\in k[[z]]$  and $\mu\in k_1[[z]]$ such that 
$$ A_1=z^{r_1}R_1(z^{m})\circ \mu^{-1}, \ \ \ \  A_2=\mu \circ z^{r_2}R_2(z^m)$$
for some integers $r_1, r_2\geq 1$ such that  $\gcd(r_1,m)=1$ and $\gcd(r_2,m)=1.$
\ec 
\pr Since the numbers $r_1,r_2$ appearing in formulas \eqref{a2}, \eqref{a3} satisfy the condition 
 $r_1r_2\equiv r\, (\mod m),$ it follows from $\gcd(r,m)=1$ that $\gcd(r_1,m)=1$ and $\gcd(r_2,m)=1.$
Moreover, since $\gcd(r_2,m)=1$ implies that \be \l{eqas} \frac{m}{\gcd(r_2,m)}=m,\ee the series $A_1$ has the required form.  
\qed 

\vskip 0.2cm


\bc \l{bc2}  Let $A\in \Gamma$ be an even formal power series. 
Then for any decomposition  $A=A_1\circ A_2$, where  $A_1,A_2\in \Gamma$, either $A_2$ is even, or there exists $\mu\in k_1[[z]]$ such that  $\mu^{-1}\circ A_2 $ is odd and $A_1\circ \mu$ is even. On the other hand, if $A$ is odd, then there exists $\mu\in k_1[[z]]$ such that $A_1\circ \mu$ and $\mu^{-1}\circ A_2$ are odd.
\ec 
\pr If $A$ is even, then $m=2$ and $r\equiv 0\ (\mod 2)$. Therefore, the condition $r_1r_2\equiv r\ (\mod m)$ implies that either $r_2\equiv 0\ (\mod 2)$, in which case $A_2$ is even, or $r_2\equiv 1\ (\mod 2)$ but $r_1\equiv 0\ (\mod 2)$,  in which case $\mu^{-1}\circ A_2 $ is odd and $A_1\circ \mu$ is even by \eqref{eqas}. On the other hand, if $A$ is odd, then $m=2$ and $r\equiv 1\ (\mod 2).$ Thus, the corollary follows from Corollary \ref{bc1}. 
\qed

\vskip 0.2cm

It was shown by Reznick in \cite{rez} that if $A\in zk[[z]]$ is a formal power series such that 
some iterate of $A$ has the form $A^{\circ s}=z^rR(z^m)$ for some $R\in zk[[z]]$ and integers $m\geq 2$, $r\geq 0$, then either $A$ itself has  a similar form, or $\ord A=1$ and $\vert A \vert$ is  finite. 
 We finish this section by showing that the part of the Reznick result concerning formal power series of order at least two is an immediate corollary of 
Theorem \ref{boee}. 

\bt \l{rez}  Let $A\in \Gamma$.  Then some iterate $A^{\circ s}$, $s\geq 1,$ has the form \linebreak $A^{\circ s}=z^rR(z^m)$ for some  $R\in k[[z]]$ and integers $m\geq 2$, $r\geq 0$ if and only if  $A=z^{r_0} R_0(z^m)$ for some $R_0\in k[[z]]$ and integer $r_0\geq 0.$ 
\et  
\pr The ``if'' part is obvious. To prove the  ``only if'' part we observe that if  $\beta_A$ is some B\"ottcher function for $A$, then  $\beta_A$  remains a   B\"ottcher function for $A^{\circ s}$, $s\geq 1.$ Thus, 
if $A^{\circ s}=z^rR(z^m)$ for some  $s\geq 1$, the ``only if'' part of Theorem \ref{boee} implies that $\beta_{A}=z L(z^m)$ for some $L\in k_0[[z]]$. Using now the ``if'' part, we conclude that $A$ has the required form. \qed

\vskip 0.2cm

Let us mention that  for every $m\geq 2$ there exist series $A\in \Gamma$ that do not have the form $\mu\circ z^rR(z^m)$ 
for some  $\mu\in k_1[[z]]$ 
but have compositional right factors of this form.  Indeed, arguing as in the proof of 
Theorem \ref{dee}, one can easily see that a composition  of  series  
$ A=A_1\circ z^{r_2}R_2(z^m)$ with $\gcd(r_2,m)=1$   
has the form $\mu \circ z^{r}R(z^m)$ for some  $\mu\in k_1[[z]]$ if and only if $A_1$ has the form $\mu\circ z^{r_{1}}R_{1}(z^m).$ Thus, if $A_1$ does  not have such a form, the same is true for $A$. 

Notice that for series $A$ as above some transition functions  have the form $z M(z^m)$ and some do not. Indeed,   all functions $\phi_A$ cannot have the form $z M(z^m)$ by Theorem \ref{traa}, but some of them have this form by Corollary \ref{cd}. Since $G_A$ is a cyclic group, this gives us examples of series of order one for which Theorem \ref{rez} is not true.

\section{Functional equation $X\circ A=Y\circ B$  and reversibility}

\subsection{Functional equation $X\circ A=Y\circ B$} 
We start this section by proving Theorem \ref{mt3} and Theorem \ref{51}. 

\vskip 0.2cm 
\noindent {\it Proof of Theorem \ref{mt3}}.  If 
\be \l{fff} X\circ A=Y\circ B,\ee 
has a solution, then setting $$F=X\circ A=Y\circ B$$ and applying Corollary \ref{iit}, we see that 
\be \l{g} \phi_A\circ\phi_B= \phi_B\circ \phi_A \ee
 for all $\phi_A\in G_A$ and  $\phi_B\in G_B.$

To prove the ``if'' part, let us observe that Lemma \ref{fif} implies that condition \eqref{g} is equivalent to the condition that  
$$ \phi _{A_{\mu}}\circ\phi _{B_{\mu}}= \phi _{B_{\mu}}\circ \phi _{A_{\mu}} $$ 
 for all $\phi _{A_{\mu}}\in G _{A_{\mu}}$ and  $\phi _{B_{\mu}}\in G _{B_{\mu}}$ for some $\mu\in k_1[[z]]$. 
 Similarly, equation \eqref{fff} has a solution for $A$ and $B$ if and only if it has a solution for $A_{\mu}$ and $B_{\mu}$ 
for some $\mu\in k_1[[z]]$. Thus,  conjugating $A$ and $B$ by $\mu=\beta_A$, without loss of generality we can assume that $A=z^n,$ $n\geq 2.$

Applying Lemma \ref{00} to equality \eqref{g} for $\phi_A= \phi_{z^n}=\v_nz$, where $\v_n\in U_n^P$, we see  that  any  $ \phi_B\in G_B$ has the form $\phi_B=z M(z^n)$ for some $M\in k_0[[z]]$. By Theorem \ref{traa}, this yields that $B$ 
 has the form $B=\mu \circ z^rR(z^n)$ for some $\mu\in k_1[[z]]$, $R\in k[[z]]$,  and  $ r\geq 0$. Therefore, equality \eqref{fff} holds for 
$$X= z^rR^n(z), \ \ \ \ Y=z^n\circ \mu^{-1}. \eqno{\Box}$$

\vskip 0.2cm
\noindent{\it Proof of Theorem \ref{51}.} Let us set $n=\ord A,$ $m=\ord B$. If \eqref{ggg} 
 has a solution in $X,Y\in zk[[z]]$ for all $s,l\geq 1$, then 
by Theorem \ref{mt3} the transition functions  
\be \l{trans} \phi_{A^{\circ l}}=\beta_A\circ \v_{nl}z\circ \beta_A^{-1}, \ \ \ \phi_{B^{\circ s}}=\beta_B\circ \v_{ms}z\circ \beta_B^{-1},\ \ \ s,l\geq 1,\ee where $\v_{nl}\in U_{nl}^P$ and $\v_{ms}\in U_{ms}^P$, commute, implying that 
$$(\beta_B^{-1}\circ \beta_A \circ \v_{nl}z \circ 
\beta_A^{-1}\circ \beta_B)  
\circ \v_{ms}z= \v_{ms}z \circ (\beta_B^{-1}\circ \beta_A \circ \v_{nl}z \circ 
\beta_A^{-1}\circ \beta_B).$$ Fixing now $l$ and $\v_{nl}$ and applying Lemma \ref{00}, we see that for every $s\geq 1$ there exists $R_s\in k[[z]]$ such that  
$$\beta_B^{-1}\circ \beta_A \circ \v_{nl}z \circ 
\beta_A^{-1}\circ \beta_B=zR_s(z^{ms}).$$
Clearly, this is possible only if 
$$\beta_B^{-1}\circ \beta_A \circ \v_{nl}z \circ 
\beta_A^{-1}\circ \beta_B=cz,$$ for some $c\in k^*$, 
 and comparing coefficients in the parts of this equality we conclude that 
$$\beta_B^{-1}\circ \beta_A \circ \v_{nl}z \circ 
\beta_A^{-1}\circ \beta_B= \v_{nl}z.$$

The last equality implies that $\beta_B^{-1}\circ \beta_A$ commutes with $\v_{nl}z$. Since this is true for every $l\geq 1$ and $\v_{nl}\in U_{nl}^P$, using again Lemma \ref{00}, we conclude that for every $l\geq 1$ there exists $M_l\in k_0[[z]]$ such that   
$$\beta_B^{-1}\circ \beta_A=zM_l(z^{nl}),$$
implying that $\beta_A=\beta_B\circ cz$ for some $c\in k^*.$  

In the other direction, it is easy to see that if $\beta_A=\beta_B\circ cz$ for some $c\in k^*,$ then for all $ s,l\geq 1$
 transition functions \eqref{trans} commute, implying by Theorem \ref{mt3} that \eqref{ggg}  has a solution. \qed

\vskip 0.2cm

Theorem \ref{51} implies the following result, obtained by Dorfer and Woracek (see \cite{ew}, Proposition 3.11).

\bc \l{wor} Let $A,B \in \Gamma$  be formal power series, and 
$\beta_A$, $\beta_B$ some B\"otcher functions. Then $A$ and $B$  commute if and only if $\beta_A=\beta_B\circ \v z$ for some $\v$ satisfying $$\v^{(\ord A-1)(\ord B-1)}=1.$$
\ec
\pr  Let us set $n=\ord A,$ $m=\ord B$. If $A$ and $B$  commute, then for all $ s,l\geq 1$ the iterates $A^{\circ l}$ and $B^{\circ s}$ also commute, implying that \eqref{ggg} has the solution $X=B^{\circ s}$, $Y=A^{\circ l}$. Thus,  $\beta_A=\beta_B\circ cz$ for some $c\in k^*$ by Theorem \ref{51}. Furthermore, since 
\be \l{sra} A=
\beta_B\circ cz\circ  z^n\circ c^{-1}z\circ \beta_B^{-1}, \ \ \ B=
\beta_B\circ z^m\circ \beta_B^{-1},\ee it follows from the commutativity of $A$ and $B$ that 
$$c^{-(n-1)}=c^{-(n-1)m}.$$ On the other hand, if  $\beta_A=\beta_B\circ c z$ for some $c$ satisfying $c^{(n-1)(m-1)}=1,$ then \eqref{sra} implies that $A$ and $B$ commute. \qed  

\subsection{Right reversibility of subsemigroups of $\Gamma$} 

Let us recall that a semigroup $S$ is called {\it right amenable} if it admits a finitely additive probability measure $\mu$ defined on all the subsets of $S$  
such that for all $a\in S$ and $T\subseteq S$ 
the equality \be \l{perd0} \mu(Ta^{-1}) = \mu(T)\ee holds, where   
 the set $Ta^{-1}$ is defined by the formula  $$ Ta^{-1}=\{s \in S\, | \,sa \in T\}.$$ 
A semigroup $S$ is called {\it right reversible}  if for all $a,b\in S$ the left ideals 
$Sa$ and  $Sb$ have a non-empty intersection, that is, if for all $a,b\in S$  there exist $x,y\in S$ such that 
$ xa=yb.$ 
It is well known and follows easily from the definition (see \cite{pater}, Proposition 1.23) that every right amenable semigroup is right reversible. 

The problems of describing right reversible and right amenable semigroups of polynomials and rational functions have been studied in the recent papers \cite{peter}, \cite{peter2}, \cite{amen}. Some analogues of the results of these papers for  finitely generated subsemigroups of $\Gamma$ were obtained in the paper \cite{amf}, mentioned in the introduction. The approach of   \cite{amf} relies on the results of \cite{sha}, for which the assumption that $S$ is    finitely generated is essential. Theorem \ref{51} provides another approach to the problem, which works equally well for infinitely generated subsemigroups of $\Gamma$.  Specifically, Theorem \ref{51} implies the following result, which contains Theorem \ref{mt4} from the introduction.

\bt \l{mt5} Every right reversible subsemigroup  $S$  of $\Gamma$ is conjugate to a subsemigroup of $\f Z.$ 
In particular, every right amenable subsemigroup  $S$  of $\Gamma$ is conjugate to a subsemigroup of $\f Z.$ 
\et 
\pr Let us fix an arbitrary element $A$ of $S$. Then for every $B\in S$ and all $s,l\geq 1$, we can apply the  right reversibility condition to the elements $A^{\circ l}$ and $ B^{\circ s}$ of $ S$ concluding that  
there exist $X,Y\in S$ such that equality \eqref{ggg} holds. Therefore, 
by Theorem \ref{51}, for every $B\in S$ the equality  $\beta_A=\beta_B\circ cz$ holds for some $c\in k^*,$ 
implying that 
$$\beta_A^{-1}\circ B\circ \beta_A= (\beta_B\circ cz)^{-1}\circ B \circ (\beta_B\circ cz)=
c^{-1}z\circ \beta_B^{-1}\circ (B\circ \beta_B)\circ cz=
$$ 
$$=c^{-1}z\circ \beta_B^{-1}\circ (\beta_B\circ z^m)\circ cz=c^{m-1}z^m,
$$
where $m=\ord B.$ 
Thus, the semigroup 
$\beta_A^{-1}\circ S\circ \beta_A$ is a subsemigroup of $\f Z.$ \qed

\end{document}